\documentclass[8pt]{article}
\usepackage{latexsym}
\usepackage{amsfonts}
\usepackage[active]{srcltx} 
\usepackage{color}
\usepackage{amsmath}
\usepackage{amssymb}
\usepackage{textcomp}

\newcommand\tstrut{\rule{0pt}{2.6ex}}

\newcommand{\ba}{\begin{array}}\newcommand{\ea}{\end{array}}
\newcommand{\sym}{S}
\newcommand{\ns}{\rm}
\renewcommand{\frak}{\mathfrak}

\newcommand{\nse}{\kern-3pt\ns$=$}\newcommand{\qd}{\hfill$\Box$\medbreak}

\newcommand{\ext}{\raise1pt\hbox{$\ts\bigwedge$}}
\renewcommand{\sym}{S}
\newcommand{\ts}{\textstyle}
\newcommand{\rf}[1]{(\ref{#1})}
\newcommand{\chii}{\raise2pt\hbox{$\chi$}}

\newcommand{\Fg}{\mbox{${\cal F}\kern-2pt_g$}}

\newcommand{\Mg}{\mbox{${\cal M}\kern-2pt_g$}}
\newcommand{\Ng}{\mbox{${\cal N}\kern-2pt_g$}}
\newcommand{\V}{V\kern-1pt}

\newcommand{\Gg}{\mbox{${\cal G}\kern-2pt_g$}}
\newcommand{\cir}{\raise1.6pt\hbox{\footnotesize$\circ$}}

\newcommand{\Res}[2]{\hbox{\ns Res}\kern-16pt\lower5pt\hbox{\footnotesize$_{#1}$}\kern2pt\left[#2\right]}
\newcommand{\qk}{quaternion-K\"ahler\kern2pt}\renewcommand{\,}{\kern1pt}

\newcommand{\End}{{\rm End}}

\renewcommand{\span}{{\rm span}}

\newcommand{\dirac}{/\kern-5pt\partial}

\newcommand{\lra}{\longrightarrow}
\renewcommand{\ts}{\textstyle}
0
\newtheorem{theo}{Theorem}[section]
\newtheorem{lemma}{Lemma}[section]

\newtheorem{prop}{Proposition}[section]\def\frac#1#2{{#1\over#2}}
\def\be#1\ee{\begin{equation}#1\end{equation}}

\allowdisplaybreaks

\begin{document}
\title{Centralizers of spin
subalgebras}
\author{
Gerardo Arizmendi\footnote{Centro de
Investigaci\'on en Matem\'aticas, A. P. 402,
Guanajuato, Gto., C.P. 36000, M\'exico. E-mail: gerardo@cimat.mx}
\footnote{Partially supported by a
CONACYT scholarship}\,\,\,\,\,\,\,
and
Rafael Herrera\footnote{Centro de
Investigaci\'on en Matem\'aticas, A. P. 402,
Guanajuato, Gto., C.P. 36000, M\'exico. E-mail: rherrera@cimat.mx}
\footnote{Partially supported by
grants from CONACyT and LAISLA (CONACyT-CNRS)
}
}

\date{\today}

\maketitle

{
\abstract{

We determine the centralizers of certain isomorphic copies of spin subalgebras $\mathfrak{spin}(r)$ in
$\mathfrak{so}(d_rm)$, where $d_r$ is the dimension of a real irreducible representation of $Cl_r^0$, the even
Clifford algebra determined by the positive definite inner product on $\mathbb{R}^r$, where
$r, m\in\mathbb{N}$.

}
}

\allowdisplaybreaks

\section{Introduction}

In this paper, we determine the centralizer subalgebras of (the isomorphic images under certain monomorphisms
of)
subalgebras $\mathfrak{spin}(r)$
in $\mathfrak{so}(d_rm)$, where $d_r$ is  the dimension of the irreducible representations of $Cl_r^0$, the
even Clifford algebra determined by $\mathbb{R}^r$ endowed with the standard positive definite inner product, 
and $r, m\in\mathbb{N}$.
The need to determine such centralizers has arisen in various geometrical settings such as the following:
\vspace{-.1in}
\begin{itemize}
 \item The holonomy algebra of Riemannian manifolds endowed with a parallel even Clifford
structure \cite{Moroianu-Semmelmann}.
 \item The automorphism group of manifolds with (almost) even Clifford (hermitian) structures
\cite{Arizmendi-Herrera-Santana}. The centralizers determined in this paper help generalize the
results on automorphisms groups of Riemannian manifolds \cite{wang,wakakuwa}, almost hermitian manifolds
\cite{tanno}, and almost
quaternion-hermitian manifolds \cite{Santana}.
 \item The structure group of Riemannian manifolds admitting twisted spin structures carrying
pure spinors  \cite{Herrera-Santana}. More precisely,
if $M$ is a smooth oriented Riemannian manifold, 
$F$ is an auxiliary Riemannian vector bundle of rank $r$,
$S(TM)$ and $S(F)$ are the locally defined spinor vector bundles of $M$ and $F$ respectively, 
$(f_1,\cdots ,f_r)$ is a local orthonormal frame of $F$,
and $m\in\mathbb{N}$ is such that the bundle $S(TM)\otimes S(F)^{\otimes m}$ is globally defined, a {\em pure
spinor field} $\phi\in\Gamma(S(TM)\otimes S(F)^{\otimes m})$ is a spinor such that its  
local $2$-forms
$\eta_{kl}^\phi(X,Y) = \left<X\wedge Y\cdot \kappa_{r*}^m(f_kf_l)\cdot \phi,\phi\right>$
induce at each point $x\in M$ a representation of $Cl_r^0$ on $T_xM$ without trivial summands. The
centralizers
determined in this paper are the orthogonal complements of $\mathfrak{spin}(r)$ in the annihilator algebra of
such a spinor. Should the spinor be parallel, such annihilator will contain the holonomy
algebra of the manifold and thus be related to the special holonomies of the Berger-Simons holonomy list
\cite{Berger,Simons}. 
\end{itemize}

The paper is organized as follows. In Section \ref{sec: prels} we recall some background material and prove  
three results which will be required later in the main theorems.
More precisely, in Subection \ref{sec:preliminaries}, we recall standard material about
Clifford algebras, Spin groups, Spin algebras, and their representations. 
In Subsection \ref{sec: calculations} we 
find explicit descriptions of the real
$\mathfrak{spin}(r)$ representations $\tilde\Delta_r$, decompositions into irreducible summands of
$\tilde\Delta_r\otimes\tilde\Delta_r$, and calculate various basic centralizers.
In Section \ref{sec: centralizers}, we prove the main results of the paper, Theorems \ref{theo: 1} and
\ref{theo: 2}.
Namely, 
in Subsection \ref{sec: centralizer r not 4}, we find the centralizers  of $\mathfrak{spin}(r)$ in
$\mathfrak{so}(d_rm)$ for $r\not \equiv 0$
$(\mbox{mod 4})$ (cf. Theorem \ref{theo: 1}) and, in Subsection \ref{sec:centralizer-Spin(4k)}, we find the
centralizers of
$\mathfrak{spin}(r)$ in $\mathfrak{so}(d_rm_1+d_rm_2)$ for $r \equiv 0$ $(\mbox{mod 4})$ (cf. Theorem
\ref{theo: 2}). The proofs
involve Riemannian homogeneous spaces, representation theory and Clifford algebras.
The separation into two cases is due to the existence of exactly one and two irreducible
representations of $Cl_r^0$ for
$r \not \equiv 0$ $(\mbox{mod } 4)$ and $r
\equiv 0$ $(\mbox{mod } 4)$ respectively.

{\em Acknowledgements}. The second named author would like to thank the International Centre for
Theoretical Physics and the Institut des Hautes \'Etudes Scientifiques for their hospitality and
support.

\section{Preliminaries}\label{sec: prels}

\subsection{Clifford algebra, spin groups and representations}\label{sec:preliminaries}

In this section we recall material that can also be consulted in \cite{Friedrich, Lawson}.
Let $Cl_n$ denote the Clifford algebra generated by all the products of the orthonormal vectors
$e_1, e_2, \ldots, e_n\in \mathbb{R}^n$ subject to the relations
\begin{eqnarray*}
e_j e_k + e_k e_j &=& -2\,\delta_{jk}, 
 \quad\mbox{for $1\leq j, k\leq n$.} 
\end{eqnarray*}
We will often write
\[e_{1\ldots s}:= e_1e_2\cdots e_s.\]
Let
\[\mathbb{C}l_n=Cl_n\otimes_{\mathbb{R}}\mathbb{C},\]
the complexification of $Cl_n$. It is well known that
\[\mathbb{C}l_n\cong \left\{
                     \begin{array}{ll}
                     \End(\mathbb{C}^{2^k}), & \mbox{if $n=2k$}\\
                     \End(\mathbb{C}^{2^k})\otimes\End(\mathbb{C}^{2^k}), & \mbox{if $n=2k+1$}
                     \end{array},
\right.
\]
where
\[\mathbb{C}^{2^k}=\mathbb{C}^2\otimes \ldots \otimes \mathbb{C}^2\]
the tensor product of $k=[{n\over 2}]$ copies of $\mathbb{C}^2$.
Let us denote
\[\Delta_n = \mathbb{C}^{2^k},\]
and consider the map
\[\kappa:\mathbb{C}l_n \lra \End(\mathbb{C}^{2^k})\]
which is an isomorphism for $n$ even and the projection onto the first summand for $n$ odd.
In order to make $\kappa_n$ explicit consider the following matrices with complex entries
\[Id = \left(\begin{array}{ll}
1 & 0\\
0 & 1
      \end{array}\right),\quad
g_1 = \left(\begin{array}{ll}
i & 0\\
0 & -i
      \end{array}\right),\quad
g_2 = \left(\begin{array}{ll}
0 & i\\
i & 0
      \end{array}\right),\quad
T = \left(\begin{array}{ll}
0 & -i\\
i & 0
      \end{array}\right).
\]
Now, consider the generators of the Clifford algebra $e_1, \ldots, e_n$ so that $\kappa_n$ can be
described as follows
\begin{eqnarray*}
e_1&\mapsto& Id\otimes Id\otimes \ldots\otimes Id\otimes Id\otimes g_1\nonumber\\
e_2&\mapsto& Id\otimes Id\otimes \ldots\otimes Id\otimes Id\otimes g_2\nonumber\\
e_3&\mapsto& Id\otimes Id\otimes \ldots\otimes Id\otimes g_1\otimes T\nonumber\\
e_4&\mapsto& Id\otimes Id\otimes \ldots\otimes Id\otimes g_2\otimes T\nonumber\\
\vdots && \dots\nonumber\\
e_{2k-1}&\mapsto& g_1\otimes T\otimes \ldots\otimes T\otimes T\otimes T\nonumber\\
e_{2k}&\mapsto& g_2\otimes T\otimes\ldots\otimes T\otimes T\otimes T,\nonumber
\end{eqnarray*}
and the last generator
\[ e_{2k+1}\mapsto i\,\, T\otimes T\otimes\ldots\otimes T\otimes T\otimes T\]
if $n=2k+1$.

Let
\[u_{+1}={1\over \sqrt{2}}(1,-i),\quad u_{-1}={1\over \sqrt{2}}(1,i)\]
which forms an orthonormal basis of $\mathbb{C}^2$ with respect to the standard Hermitian product.
Note that
\[g_1(u_{\pm1})= iu_{\mp1},\quad
g_2(u_{\pm1})= \pm u_{\mp1},\quad
T(u_{\pm1})= \mp u_{\pm1}.\]
Thus, we get a unitary basis of $\Delta_n=\mathbb{C}^{2^k}$
\[\mathcal{B}=\{u_{\varepsilon_1,\ldots,\varepsilon_k}=u_{\varepsilon_1}\otimes\ldots\otimes
u_{\varepsilon_k}\,\,|\,\, \varepsilon_j=\pm 1,
j=1,\ldots,k\},\]
with respect to the induced Hermitian product on $\mathbb{C}^{2^k}$.

The Clifford multiplication of a vector $e$ and a spinor $\psi$ is defined by
$ e\cdot \psi = \kappa_n(e)(\psi).$
Thus, if $1\leq j\leq k$
\begin{eqnarray*}
e_{2j-1}\cdot u_{\varepsilon_1,\ldots,\varepsilon_k}&=& i(-1)^{j-1}
\left(\prod_{\alpha=k-j+2}^k \varepsilon_{\alpha}\right)
u_{\varepsilon_1,\ldots, (-\varepsilon_{k-j+1}) ,\ldots,\varepsilon_k}    \nonumber\\
e_{2j}\cdot u_{\varepsilon_1,\ldots,\varepsilon_k}&=&  (-1)^{j-1}
\left(\prod_{\alpha=k-j+1}^k \varepsilon_{\alpha}\right)
u_{\varepsilon_1,\ldots, (-\varepsilon_{k-j+1}) ,\ldots,\varepsilon_k}   \nonumber
\end{eqnarray*}
and
\[e_{2k+1}\cdot u_{\varepsilon_1,\ldots,\varepsilon_k}= i(-1)^k
\left(\prod_{\alpha=1}^k \varepsilon_{\alpha}\right) u_{\varepsilon_1,\ldots,\varepsilon_k}\]
if $n=2k+1$ is odd.

The Spin group $Spin(n)\subset Cl_n$ is the subset
\[Spin(n) =\{x_1x_2\cdots x_{2l-1}x_{2l}\,\,|\,\,x_j\in\mathbb{R}^n, \,\,
|x_j|=1,\,\,l\in\mathbb{N}\},\]
endowed with the product of the Clifford algebra.
The Lie algebra of $Spin(n)$ is
\[\mathfrak{spin}(n)=\mbox{span}\{e_ie_j\,\,|\,\,1\leq i< j \leq n\}.\]
The restriction of $\kappa$ to $Spin(n)$ defines the Lie group representation
\[\kappa_n:=\kappa|_{Spin(n)}:Spin(n)\lra GL(\Delta_n),\]
which is, in fact, special unitary \cite{Friedrich}.

There exist either real or quaternionic structures on the spin representations.
A quaternionic structure $\alpha$ on $\mathbb{C}^2$ is given by
\[\alpha\left(\begin{array}{c}
z_1\\
z_2
              \end{array}
\right) = \left(\begin{array}{c}
-\overline{z}_2\\
\overline{z}_1
              \end{array}\right),\]
and a real structure $\beta$ on $\mathbb{C}^2$ is given by
\[\beta\left(\begin{array}{c}
z_1\\
z_2
              \end{array}
\right) = \left(\begin{array}{c}
\overline{z}_1\\
\overline{z}_2
              \end{array}\right).\]
Note that these structures satisfy
\[
\begin{array}{rclcrcl}
 \left< \alpha(v),w\right> &=& \overline{\left< v,\alpha(w)\right> }, &\quad&
 \left< \alpha(v),\alpha(w)\right> &=& \overline{\left< v,w\right> }, \\
 \left< \beta(v),w\right> &=& \overline{\left< v,\beta(w)\right> }, &\quad&
 \left< \beta(v),\beta(w)\right> &=& \overline{\left< v,w\right> }, 
\end{array}
\]
with respect to the standard hermitian product in $\mathbb{C}^2$, where $v,w\in \mathbb{C}^2$.
The real and quaternionic structures $\gamma_n$  on $\Delta_n=(\mathbb{C}^2)^{\otimes
[n/2]}$ are built as follows
\[
\begin{array}{cclll}
 \gamma_n &=& (\alpha\otimes\beta)^{\otimes 2k} &\mbox{if $n=8k,8k+1$}& \mbox{(real),} \\
 \gamma_n &=& \alpha\otimes(\beta\otimes\alpha)^{\otimes 2k} &\mbox{if $n=8k+2,8k+3$}&
\mbox{(quaternionic),} \\
 \gamma_n &=& (\alpha\otimes\beta)^{\otimes 2k+1} &\mbox{if $n=8k+4,8k+5$}&\mbox{(quaternionic),} \\
 \gamma_n &=& \alpha\otimes(\beta\otimes\alpha)^{\otimes 2k+1} &\mbox{if $n=8k+6,8k+7$}&\mbox{(real).}
\end{array}
\]
which also satisfy
\[
\begin{array}{rclcrcl}
 \left< \gamma_n(v),w\right> &=& \overline{\left< v,\gamma_n(w)\right> }, &\quad&
 \left< \gamma_n(v),\gamma_n(w)\right> &=& \overline{\left< v,w\right> }, \\
\end{array}
\]
where $v,w\in\Delta_n$.
This means 
\begin{equation}
 \left< v+ \gamma_n(v),w+ \gamma_n(w)\right> \in \mathbb{R}.   \label{eq: real inner product}
\end{equation}

\begin{lemma} \label{Clifford commute}
 Let $m\geq r$ and let $e_{i_1}\dots e_{i_r}=:e_{i_1 \dots i_r}=e_I\in Cl_m$. Then $e_I$ commutes with
$\mathfrak{spin}(r)=span\{e_ie_j|1\leq i<j\leq r\}$ if and only if $I\subset\{r+1,\dots,m\}$ or
$\{1,\dots,r\}\subset I$.
\end{lemma}
{\em Proof}. Suppose that neither $I\subset\{r+1,\dots,m\}$ nor $\{1,\dots,r\}\subset I$ then there exist
$j,k\in \{1,\dots,r\}$ such that $j\in I$ and $k\not \in I$. Rearranging the other of the ${i_l}$'s if
necessary we can suppose that $j=i_1$, so that $e_Ie_je_k=e_{i_1}\dots e_{i_{r}}e_{i_1}e_k=(-1)^re_{i_2}\dots 
e_{i_{r}}. e_k$ and $e_je_ke_I=e_{i_1}e_ke_{i_1}\dots e_{i_{r}}=(-1)^{r+1}e_{i_2}\dots  e_{i_{r}} e_k$.

Conversely, the volume form on $Cl_r$ commutes with $\mathfrak{spin}(r)$ in every dimension and if $k\not \in
\{1,\dots,r\}$ then for all $i,j\in \{1,\dots,r\}$  we have that $e_ie_je_k=e_ke_ie_j$.
\qd

Now, we summarize some results about real representations of $Cl_r^0$ in the next table (cf. \cite{Lawson}).
Here $d_r$ denotes the dimension of an irreducible representation of $Cl^0_r$ and $v_r$ the number of distinct
irreducible representations.
\[\begin{array}{|c|c|c|c|c|}
\hline
r \mbox{ (mod 8)}&Cl_r^0&d_r&v_r \tstrut\\
\hline
1&\mathbb R(d_r)&2^{\lfloor\frac r2\rfloor}&1 \tstrut\\
\hline
2&\mathbb C(d_r/2)&2^{\frac r2}&1 \tstrut\\
\hline
3&\mathbb H(d_r/4)&2^{\lfloor\frac r2\rfloor+1}&1 \tstrut\\
\hline
4&\mathbb H(d_r/4)\oplus \mathbb H(d_r/4)&2^{\frac r2}&2 \tstrut\\
\hline
5&\mathbb H(d_r/4)&2^{\lfloor\frac r2\rfloor+1}&1 \tstrut\\
\hline
6&\mathbb C(d_r/2)&2^{\frac r2}&1 \tstrut\\
\hline
7&\mathbb R(d_r)&2^{\lfloor\frac r2\rfloor}&1 \tstrut\\
\hline
8&\mathbb R(d_r)\oplus \mathbb R(d_r)&2^{\frac r2-1}&2 \tstrut\\
\hline
\end{array}
\]
\centerline{Table 1}

Let $\tilde\Delta_r$ denote the irreducible representation of $Cl_r^0$ for $r\not\equiv0$ $(\mbox{mod } 4) $
and $\tilde\Delta^{\pm}_r$ denote the irreducible representations for $r\equiv0$ $(\mbox{mod } 4)$. Note that
the representations are complex for $r\equiv 2,6$ $(\mbox{mod } 8)$ and quaternionic for $r\equiv 3,4,5$ 
$(\mbox{mod } 8)$. It is interesting to note that these features are reflected in the main results of the
paper.

Note also that if $r\equiv 4,6,7,8$ $(\mbox{mod } 8)$ then $d_r=d_{r-1}$ and if $r\equiv 1,2,3,5$ $(\mbox{mod
} 8)$ then $d_r=2d_{r-1}$. By restricting to a standard subagebra $Cl^0_{r-1}\subset Cl_r^0$, the
representations decompose as follows:
\[\begin{array}{|c|c|}
\hline 
r \mbox{ (mod 8)}&\tilde\Delta_r|_{Cl^0_{r-1}} \tstrut\\ 
\hline 
1&\tilde\Delta_r\cong\tilde\Delta^+_{r-1}+\tilde\Delta^-_{r-1} \tstrut\\
\hline
2&\tilde\Delta_r\cong\tilde\Delta_{r-1}+\tilde\Delta_{r-1} \tstrut\\
\hline
3&\tilde\Delta_r\cong\tilde\Delta_{r-1}+\tilde\Delta_{r-1} \tstrut\\
\hline
4&\tilde\Delta^{\pm}_r\cong\tilde\Delta_{r-1} \tstrut\\
\hline
5&\tilde\Delta_r\cong\tilde\Delta_{r-1}^++\tilde\Delta^-_{r-1} \tstrut\\
\hline
6&\tilde\Delta_r\cong\tilde\Delta_{r-1} \tstrut\\
\hline
7&\tilde\Delta_r\cong\tilde\Delta_{r-1} \tstrut\\
\hline
8&\tilde\Delta^{\pm}_r\cong\tilde\Delta_{r-1} \tstrut\\
\hline
\end{array}
\]
\centerline{Table 2}

\subsection{Real spin representations and basic centralizers}\label{sec: calculations}

In this section we prove results which are essential in Theorems \ref{theo: 1} and
\ref{theo: 2}. 
Let $\ext^2V$ and $\sym^2V$ denote the second exterior and symmetric power of a finite dimensional vector
space respectively.
In addition, if the vector space is endowed with an inner product, let $\sym_0^2 V$ denote the orthogonal
complement of the identity endomorphism within the symmetric endomorphisms of $V$.

\begin{prop}\label{prop: 1} The centralizers of the spin subalgebras under consideration are:
\[  \begin{array}{|c|c|c|}
 \hline
  r \mbox{ {\rm (mod 8)}}
&C_{\mathfrak{so}(d_r)}({\mathfrak{spin}({r})})&C_{\mathfrak{so}(d_r)\oplus\mathfrak{so}(d_r)}({\mathfrak{spin
}({r})})\\
 \hline
 0 &  &\{0\}\\
 \hline
\pm 1 & \{0\} &\\
\hline
\pm 2 & \mathfrak{u}(1) &\\
\hline
\pm 3 & \mathfrak{sp}(1) &\\
\hline
 4 &  & \mathfrak{sp}(1)\oplus\mathfrak{sp}(1)\\
\hline
  \end{array}
\]
Furthermore, the representations $\ext^2 \tilde{\Delta}_{r}$, $\ext^2 \tilde{\Delta}_{r}^{\pm}$,
$\sym_0^2\tilde{\Delta}_{r}$, $\sym_0^2\tilde{\Delta}_{r}^{\pm}$ and
$\tilde{\Delta}_{r}^+\otimes\tilde{\Delta}_{r}^-$ have the following trivial $Spin(r)$ subrepresentations:
\[
\begin{array}{|c|c|c|c|}
\hline
  r \mbox{ {\rm (mod 8)}} & \ext^2 \tilde{\Delta}_{r}, \ext^2 \tilde{\Delta}_{r}^\pm  &
\sym_0^2 \tilde{\Delta}_{r},\sym_0^2 \tilde{\Delta}_{r}^\pm &\tilde{\Delta}_{r}^+\otimes\tilde{\Delta}_{r}^-
\rule{0pt}{3ex}\\
\hline
0 & \{0\} & \{0\} & \{0\} \\
\hline
\pm 1 &  \{0\} & \{0\} &\\
\hline
\pm 2 &  \mathfrak{u}(1) & \{0\} &\\
\hline
\pm3 &\mathfrak{sp}(1) &\{0\} &\\
\hline
4 & \mathfrak{sp}(1) &\{0\} &\{0\} \\
\hline
\end{array}
\]
\end{prop}

\subsubsection*{Proof of the Proposition}

\subsubsection*{Case $r\equiv \pm1 \,\,\,({\rm mod}\,\,\, 8)$ }
In both cases there exist real structures $\gamma=\gamma_{r}$ on $\Delta_r$. 
By using these real structures, we can describe the underlying real space $\tilde\Delta_{r}\subset\Delta_{r}$
as follows.
Recall the unitary basis $\mathcal{B}$ of $\Delta_{r}$   and let
\[\mathcal{B}_1=\{u_{\varepsilon_1,\ldots,\varepsilon_{[r/2]}}+\gamma(u_{\varepsilon_1,\ldots,\varepsilon_{[
r/2 ] } } ) ,
iu_{
\varepsilon_1,\ldots,\varepsilon_{[r/2]}}+\gamma(iu_{\varepsilon_1,\ldots,\varepsilon_{[r/2]}})
\,\,|\,\, \varepsilon_j=\pm 1,\,\, j=1,\ldots,{[r/2]}\},\]
which is an orthogonal basis for
\[\tilde\Delta_{r}=\span(\mathcal{B}_1) = \{v+\gamma_r(v)\,|\, v\in\Delta_r\}.\]
since the hermitian product of $\Delta_r$
restricts to a real inner product on $\tilde\Delta_r$ (cf. \rf{eq: real inner product}).
Consider the $\mathfrak{spin}(r)$ equivariant morphism
\[\Phi:\tilde\Delta_{r}\otimes\tilde\Delta_{r}\rightarrow\bigoplus_k\ext^{2k}\mathbb{R}^{r}.\]
defined by
\[\Phi[(v+\gamma(v))\otimes (w+\gamma(w))]=\sum_{k=0}^{[r/2]}\sum_{j_1<\dots<j_{2k}} \left<e_{j_1}\dots
e_{j_{2k}}(v+\gamma(v)),(w+\gamma(w))\right>e_{j_1}\dots e_{j_{2k}},\]
where
$\left<e_{j_1}\dots e_{j_{2k}}(v+\gamma(v)),w+\gamma(w)\right>$ 
is real.
Now, let $v+\gamma(v)\in \mathcal{B}_1$ and let $\tilde{v}=e_{j_1}\dots e_{j_{2k}}v$, so that
$\pm(\tilde{v}+\gamma(\tilde{v}))\in \mathcal{B}_1$ and
\begin{eqnarray*}
\left<e_{j_1}\dots e_{j_{2k}}(v+\gamma(v)),\tilde{v}+\gamma(\tilde{v})\right>
 &=&\left<e_{j_1}\dots e_{j_{2k}}v+\gamma(e_{j_1}\dots e_{j_{2k}}v),\tilde{v}+\gamma(\tilde{v})\right>\\
 &=&\left<e_{j_1}\dots e_{j_{2k}}v+\gamma(e_{j_1}\dots e_{j_{2k}}v),e_{j_1}\dots
e_{j_{2k}}v+\gamma(e_{j_1}\dots e_{j_{2k}}v)\right>\\
 &=&2.
\end{eqnarray*}
Hence, the image $\Phi(\tilde\Delta_{r}\otimes\tilde\Delta_{r})$ has non-trivial projection to
$\ext^{2k}\mathbb{R}^{r}$ for $k=0,\ldots,{[r/2]}$. Since the dimensions of
$\tilde\Delta_{r}\otimes\tilde\Delta_{r}$ and $\bigoplus_k\ext^{2k}\mathbb{R}^{r}$ coincide, $\Phi$ is
equivariant and $\bigoplus_k\ext^{2k}\mathbb{R}^{r}$ is a sum of nonequivalent irreducible representations,
Schur's Lemma implies that
\[\tilde\Delta_{r}\otimes\tilde\Delta_{r}\cong\bigoplus_k\ext^{2k}\mathbb{R}^{r}\] as 
$\mathfrak{spin}(r)$ representations.
Moreover, since
\begin{eqnarray*}
\left<e_{j_1}\dots e_{j_{4l}}(v+\gamma(v)),w+\gamma(w)\right>
 &=&\left<v+\gamma(v),e_{j_1}\dots e_{j_{4l}}(w+\gamma(w))\right>\\
 &=&\left<e_{j_1}\dots e_{j_{4l}}(w+\gamma(w)),v+\gamma(v)\right>,
\end{eqnarray*}
$\Phi(\ext^2\tilde\Delta_{r})$ has trivial projection $\ext^{4l}\mathbb{R}^{r}$.

Analogously, since
\begin{eqnarray*}
 \left<e_{j_1}\dots e_{j_{4l+2}}(v+\gamma(v)),w+\gamma(w)\right>
 &=&\left<v+\gamma(v),-e_{j_1}\dots e_{j_{4l+2}}(w+\gamma(w))\right>\\
 &=&\left<-e_{j_1}\dots e_{j_{4l+2}}(w+\gamma(w)),v+\gamma(v)\right>,
\end{eqnarray*}
$\Phi(\sym^2\tilde\Delta_{r})$ has trivial projection to $\ext^{4l+2}\mathbb{R}^{r}$.

Recall that we wish to find the centralizer of $\mathfrak{spin}(r)$ in
$\mathfrak{so}(d_r)\subset \End(\tilde\Delta_r) \cong Cl_r^0$.
Notice that, by Lemma \ref{Clifford commute},
\[C_{{\rm End}(\tilde\Delta_r)}({\mathfrak{spin}(r)})=C_{Cl_r^0}({\mathfrak{spin}(r)})=\span({\rm
Id}_{d_r\times d_r}).\]
Hence,
\[C_{\mathfrak{so}(d_r)}({\mathfrak{spin}(r)})=\{0\}\subset\mathfrak{so}(d_r)\cong\ext^2\tilde\Delta_r.\]

\subsubsection*{Case $r\equiv -2 \,\,\,({\rm mod}\,\,\, 8)$ }

Recall that for $r\equiv -2$ $\mbox{\rm (mod 8)}$ we have 
\[\tilde\Delta_r\cong \tilde\Delta_{r+1}\]
as representations of $\mathfrak{spin}(r)$, and
\[\tilde\Delta_{r}=\span(\mathcal{B}_1).\]
\[\tilde\Delta_r\otimes \tilde\Delta_r\cong\tilde\Delta_{r+1}\otimes \tilde\Delta_{r+1}.\]
Since $r+1\equiv -1$ (mod 8),
\[\tilde\Delta_{r+1}\otimes \tilde\Delta_{r+1} \cong \ext^{ev}\mathbb{R}^{r+1}\]
with respect to $\mathfrak{spin}(r+1)$,
as proved in the previous subsection.
Furthermore,  $\mathbb{R}^{r+1} = \mathbb{R}^r \oplus 1$,
\begin{eqnarray}
 \ext^0 \mathbb{R}^{r+1} &=& 1, \nonumber\\
 \ext^2 \mathbb{R}^{r+1} &=& \ext^2 \mathbb{R}^r +  \mathbb{R}^r,\nonumber\\
 \ext^4 \mathbb{R}^{r+1} &=& \ext^4 \mathbb{R}^r + \ext^3\mathbb{R}^r, \nonumber
\\
 \vdots \nonumber\\
 \ext^{r} \mathbb{R}^{r+1} &=& \ext^{r} \mathbb{R}^r + \ext^{r-1}\mathbb{R}^r.\nonumber
\end{eqnarray}
and
\begin{eqnarray*}
\tilde\Delta_r\otimes \tilde\Delta_r
  &=&\ext^*\mathbb{R}^r.
\end{eqnarray*}
On the other hand,
\begin{eqnarray*}
\tilde\Delta_r\otimes \tilde\Delta_r
  &=& \ext^2\tilde\Delta_r + \sym_0^2 \tilde\Delta_r + 1.
\end{eqnarray*}
and
\begin{eqnarray*}
\ext^2\tilde\Delta_{r}&\cong& \ext^2\tilde\Delta_{r+1}\\
&\cong& \bigoplus\ext^{4l+2}\mathbb{R}^{r+1}, \\
&\cong& \bigoplus_{l\geq 0}\ext^{4l+1}\mathbb{R}^{r}\bigoplus_{l\geq 0}\ext^{4l+2}\mathbb{R}^{r}, \\
\sym_0^2\tilde\Delta_{r}&\cong&\sym_0^2\tilde\Delta_{r+1}\\
&\cong& \bigoplus_{l> 0}\ext^{4l}\mathbb{R}^{r+1}\\
&\cong& \bigoplus_{l> 0}\ext^{4l}\mathbb{R}^{r}\bigoplus_{l\geq 0}\ext^{4l+3}\mathbb{R}^{r}.
\end{eqnarray*}
We see that $\ext^2\tilde\Delta_r$ contains a 1-dimensional trivial $\mathfrak{spin}(r)$ representation.

Recall that we wish to find the centralizer of $\mathfrak{spin}(r)$ in
$\mathfrak{so}(d_r)\subset \End(\tilde\Delta_r) 
\cong Cl_r$.
Note that any element of $Cl_r$ which commutes with $\mathfrak{spin}(r)$ must commute with
the volume element $e_1e_2\cdots e_r \in Cl_r$, and such elements are precisely $Cl_r^0$. Thus, by Lemma
\ref{Clifford commute}
\[C_{Cl_r}({\mathfrak{spin}(r)})\subseteq C_{Cl_r^0}({\mathfrak{spin}(r)}) =
\span(1)\oplus \span( e_1\cdots e_r), \]
where
$e_1\cdots e_r$ acts as an orthogonal complex structure $J$ on $\tilde\Delta_r$ 
which generates the afore mentioned 1-dimensional trivial summand in $\ext^2\tilde\Delta_r$. 
Hence, 
\[C_{\mathfrak{so}(d_r)}(\mathfrak{spin}(r)) = {\rm span}(J) \cong 
\mathfrak{u}(1)\subset\mathfrak{so}(d_r)\cong\ext^2\tilde\Delta_r.\]

\subsubsection*{Case $r\equiv 2 \,\,\,({\rm mod}\,\,\, 8)$ }

In this case, there exists a quaternionic structure $\gamma_{r}$ on
$\Delta_r$which commutes with Clifford
multiplication. We can describe the real space
$\tilde\Delta_{r}\subset\Delta_{r}$ as follows.
Recall the unitary basis $\mathcal{B}$ of $\Delta_{r}$  and let
\[\mathcal{B}_2=\{u_{\varepsilon_1,\ldots,\varepsilon_{r/2}}+\gamma(u_{\varepsilon_1,\ldots,\varepsilon_{r/2}}
) ,
iu_{
\varepsilon_1,\ldots,\varepsilon_{r/2}}+\gamma(iu_{\varepsilon_1,\ldots,\varepsilon_{r/2}})
\,\,|\,\, \varepsilon_j=\pm 1,\,\, j=1,\ldots,{r/2}\}.\]
Note that the space generated by the orthogonal basis $\mathcal{B}_2$ is preserved by the action of
$\mathfrak{spin}(r)$ and $Cl_r^0$, the hermitian product in $\Delta_r$ restricts to an inner product to
$\tilde\Delta_r$ (cf. \rf{eq: real inner product}), and its dimension is $d_r$. Therefore 
\[\tilde\Delta_{r}=\span(\mathcal{B}_2).\] 
Now consider the $\mathfrak{spin}(r)$ equivariant morphism
\[\Phi:\tilde\Delta_{r}\otimes\tilde\Delta_{r}\rightarrow\bigoplus_k\ext^{k}\mathbb{R}^{r}\]
defined by
\[\Phi[(v+\gamma(v))\otimes(w+\gamma(w))]=\sum_{k=0}^{r}\sum_{j_1<\dots<j_{k}} \left<e_{j_1}\dots
e_{j_{2k}}(v+\gamma(v)),(w+\gamma(w))\right>e_{j_1}\dots e_{j_{k}},\]
where
$\left<e_{j_1}\dots
e_{j_{2k}}(v+\gamma(v)),w+\gamma(w)\right>$ is real.
Let $v+\gamma(v)\in \mathcal{B}_2$ and $\tilde{v}=e_{j_1}\dots e_{j_{k}}v$, then
$\pm(\tilde{v}+\gamma(\tilde{v}))\in \mathcal{B}_2$ and
\begin{eqnarray*}
\left<e_{j_1}\dots e_{j_{k}}(v+\gamma(v)),\tilde{v}+\gamma(\tilde{v})\right>
 &=&\left<e_{j_1}\dots e_{j_{k}}v+\gamma(e_{j_1}\dots e_{j_{k}}v),\tilde{v}+\gamma(\tilde{v})\right>\\
 &=&\left<e_{j_1}\dots e_{j_{k}}v+\gamma(e_{j_1}\dots e_{j_{k}}v),e_{j_1}\dots e_{j_{k}}v+\gamma(e_{j_1}\dots
e_{j_{k}}v)\right>\\
 &=&2.
\end{eqnarray*}
Hence, the image $\Phi(\tilde\Delta_{r}\otimes\tilde\Delta_{r})$ has non-trivial projection to
$\ext^{k}\mathbb{R}^{r}$ for $k=0,\ldots,{r}$. Since the dimensions of
$\tilde\Delta_{r}\otimes\tilde\Delta_{r}$ and $\bigoplus_k\ext^{k}\mathbb{R}^{r}$ coincide, $\Phi$ is
equivariant and $\bigoplus_k\ext^{k}\mathbb{R}^{r}$ is a sum of nonequivalent irreducible representations,
Schur's Lemma implies that
\[\tilde\Delta_{r}\otimes\tilde\Delta_{r}\cong\bigoplus_k\ext^{k}\mathbb{R}^{r}\] 
as $\mathfrak{spin}(r)$ representations.
Moreover,
\begin{eqnarray*}
\left<e_{j_1}\dots e_{j_{4l}}(v+\gamma(v)),w+\gamma(w)\right>
 &=&\left<v+\gamma(v),e_{j_1}\dots e_{j_{4l}}(w+\gamma(w))\right>\\
 &=&\left<e_{j_1}\dots e_{j_{4l}}(w+\gamma(w)),v+\gamma(v)\right>,
\end{eqnarray*}
and
\begin{eqnarray*}
\left<e_{j_1}\dots e_{j_{4l+3}}(v+\gamma(v)),w+\gamma(w)\right>
 &=&\left<v+\gamma(v),e_{j_1}\dots e_{j_{4l+3}}(w+\gamma(w))\right>\\
 &=&\left<e_{j_1}\dots e_{j_{4l+3}}(w+\gamma(w)),v+\gamma(v)\right>.
\end{eqnarray*}
Therefore $\Phi(\ext^2\tilde\Delta_{r})$ has trivial projection to $\ext^{4l}\mathbb{R}^{r}$ and
$\ext^{4l+3}\mathbb{R}^{r}$.
Analogously,
\begin{eqnarray*}
 \left<e_{j_1}\dots e_{j_{4l+1}}(v+\gamma(v)),w+\gamma(w)\right>
 &=&\left<v+\gamma(v),-e_{j_1}\dots e_{j_{4l+1}}(w+\gamma(w))\right>\\
 &=&\left<-e_{j_1}\dots e_{j_{4l+1}}(w+\gamma(w)),v+\gamma(v)\right>,
\end{eqnarray*}
and
\begin{eqnarray*}
 \left<e_{j_1}\dots e_{j_{4l+2}}(v+\gamma(v)),w+\gamma(w)\right>
 &=&\left<v+\gamma(v),-e_{j_1}\dots e_{j_{4l+2}}(w+\gamma(w))\right>\\
 &=&\left<-e_{j_1}\dots e_{j_{4l+2}}(w+\gamma(w)),v+\gamma(v)\right>.
\end{eqnarray*}
Therefore,
$\Phi(\sym^2\tilde\Delta_{r})$ has trivial projection to $\ext^{4l+1}\mathbb{R}^{r}$ and
$\ext^{4l+2}\mathbb{R}^{r}$.

Recall that we wish to compute the centralizer of $\mathfrak{spin}(r)$ in 
$\mathfrak{so}(d_r)\subset \End(\tilde\Delta_r)\cong Cl_r$.
As in the previous case,
any element of $Cl_r$ which commutes with $\mathfrak{spin}(r)$ must commute with
the volume element $e_1e_2\cdots e_r \in Cl_r$, and such elements are precisely $Cl_r^0$. Thus, by
Lemma \ref{Clifford commute},
\[C_{Cl_r}({\mathfrak{spin}(r)})\subseteq C_{Cl_r^0}({\mathfrak{spin}(r)}) = \span(1)\oplus \span( e_1\cdots
e_r).\]
where
$e_1\cdots e_r$ acts as an orthogonal complex structure $J$ on $\tilde\Delta_r$. Hence, 
\[C_{\mathfrak{so}(d_r)}(\mathfrak{spin}(r)) = {\rm span}(J) \cong 
\mathfrak{u}(1)\subset\mathfrak{so}(d_r)\cong\ext^2\tilde\Delta_r.\]

\subsubsection*{Case $r\equiv 3 \,\,\,({\rm mod}\,\,\, 8)$ }

Recall that \[\tilde\Delta_r\oplus \tilde\Delta_r\cong
\tilde\Delta_{r+1}^+\oplus\tilde\Delta_{r+1}^-\cong\tilde\Delta_{r+2}\cong\tilde\Delta_{r+3}\]
as representations of $\mathfrak{spin}(r)$.
Since $\gamma_{r+3}$ is
a real structure, 
\[\tilde\Delta_{r+3}=\{v+\gamma_{r+3}(v) \,\,|\,\, v\in\Delta_{r+3}\}.\]
Moreover, 
\[\tilde\Delta_r={1\over 2}(1\pm e_1\dots e_r)\tilde\Delta_{r+3},\]
so that
\begin{eqnarray*}
\tilde\Delta_{r+3}\otimes \tilde\Delta_{r+3}&\cong&(\tilde\Delta_r\oplus \tilde\Delta_r)\otimes
(\tilde\Delta_r\oplus \tilde\Delta_r)\\
&\cong&4\,\,\tilde\Delta_r\otimes \tilde\Delta_r
\end{eqnarray*}
with respect to $\mathfrak{spin}(r)$.
Since $r+3\equiv -2$ (mod 8),
\[\tilde\Delta_{r+3}\otimes \tilde\Delta_{r+3} \cong \ext^{*}\mathbb{R}^{r+3}\]
with respect to $\mathfrak{spin}(r+3)$.
Now $\mathbb{R}^{r+3} = \mathbb{R}^r \oplus 3$,
\begin{eqnarray*}
 \ext^0 \mathbb{R}^{r+3} &=& 1,\\
 \ext^1 \mathbb{R}^{r+3} &=&  \mathbb{R}^r + 3,\\
 \ext^2 \mathbb{R}^{r+3} &=& \ext^2 \mathbb{R}^r +  3\mathbb{R}^r+3,\\
 \ext^3 \mathbb{R}^{r+3} &=& \ext^3 \mathbb{R}^r + 3\ext^2\mathbb{R}^r+3\mathbb{R}^r+1,\\
 \vdots\\
 \ext^{r+3} \mathbb{R}^{r+3} &=& \ext^{r+3} \mathbb{R}^r + 3\ext^{r+2}\mathbb{R}^r+3\ext^{r+1}\mathbb{R}^r+
\ext^{r}\mathbb{R}^r=1,
\end{eqnarray*}
and we have 
\begin{eqnarray*}
4\tilde\Delta_r\otimes \tilde\Delta_r
  &=& 8\ext^*\mathbb{R}^r.
\end{eqnarray*}
Therefore \begin{eqnarray*}
\tilde\Delta_r\otimes \tilde\Delta_r
  &=& 2\ext^*\mathbb{R}^r.
\end{eqnarray*}

Recall that we wish to compute the centralizer of $\mathfrak{spin}(r)$ in 
$\mathfrak{so}(d_r)\subset \End(\tilde\Delta_r)\subset\End(\tilde\Delta_{r+3})\cong Cl_{r+3}$.
First, we will compute $C_{Cl_{r+3}}(\mathfrak{spin}(r))$.
Suppose
\[\eta = \sum_{I} \eta_Ie_I \in C_{Cl_{r+3}}(\mathfrak{spin}(r)),\]
it must commute in Clifford product with  every $e_ie_j\in \mathfrak{spin}(r)$, $1\leq i<j\leq r$.
By Lemma \ref{Clifford commute}, 
the only free coefficients are
$\eta_{\emptyset}$,
$\eta_{r+1}$,
$\eta_{r+2}$,
$\eta_{r+3}$,
$\eta_{r+1,r+2}$,
$\eta_{r+1,r+3}$,
$\eta_{r+2,r+3}$,
$\eta_{r+1,r+2,r+3}$,
$\eta_{1,\dots,r}$,
$\eta_{1,\dots,r+1}$,
$\eta_{1,\dots,r,r+2}$,
$\eta_{1,\dots,r,r+3}$,
$\eta_{1,\dots,r+2}$,
$\eta_{1,\dots,r+1,r+3}$,
$\eta_{1,\dots,r,r+2,r+3}$,
$\eta_{1,\dots,r+3}$,
i.e.
\begin{eqnarray*}
C_{Cl_{r+3}}(\mathfrak{spin}(r)) &=& \span(1,
e_{r+1},
e_{r+2},
e_{r+3},
e_{r+1,r+2},
e_{r+1,r+3},
e_{r+2,r+3},
e_{r+1,r+2,r+3},
e_{1,\dots,r},\\
&& 
e_{1,\dots,r+1},
e_{1,\dots,r,r+2},
e_{1,\dots,r,r+3},
e_{1,\dots,r+2},
e_{1,\dots,r+1,r+3},
e_{1,\dots,r,r+2,r+3},
e_{1,\dots,r+3}.
)
\end{eqnarray*}
These elements act as automorphisms of $\tilde\Delta_{r+3}$ via Clifford multiplication.
Recall that the two copies of (the $d_r$-dimensional real representation) $\tilde\Delta_{r}$ are
\begin{eqnarray*}
\tilde\Delta_{r} &=& {1\over 2} (1+e_1\dots e_{r})\cdot \tilde\Delta_{r+3},\\
\tilde\Delta_{r} &=& {1\over 2} (1-e_1\dots e_{r})\cdot \tilde\Delta_{r+3}.
\end{eqnarray*}
We will restrict our attention to the first copy. In order to project these elements to automorphisms of this
copy we need to multiply by ${1\over 2} (1+e_1\dots e_{r})$. Observe for example that
\begin{eqnarray*}
{1\over 2} (1+e_1\dots e_{r})\cdot 1 &=&{1\over 2} (1+e_1\dots e_{r}),\\
{1\over 2} (1+e_1\dots e_{r})\cdot e_1\dots e_{r}&=&{1\over 2} (1+e_1\dots e_{r}),\\
{1\over 2} (1+e_1\dots e_{r})\cdot  e_{r+1}&=&{1\over2}(e_{r+1}+e_{1,\dots,r+1}),\\
{1\over 2} (1+e_1\dots e_{r})\cdot   e_1\dots e_{r+1}&=&{1\over2}(e_{r+1}+e_{1,\dots,r+1}),
\end{eqnarray*}
and, for $v\in\Delta_{r+3}$,
\begin{eqnarray*}
{1\over 2}
(1+e_1\dots e_r)\cdot {1\over 2}
(1+e_1\dots e_r)\cdot(v+\gamma_{r+3}(v))&=&{1\over 2}
(1+e_1\dots e_r)\cdot(v+\gamma_{r+3}(v))\\
{1\over 2}
((e_{r+1}+e_{1,\dots,r+1}))\cdot {1\over 2}
(1+e_1\dots e_r)\cdot(v+\gamma_{r+3}(v))&=&0, 
\end{eqnarray*}
so 
${1\over 2} (1+e_1\dots e_{r})$ acts as the identity element on this copy of $\tilde\Delta_r$ and
${1\over2}(e_{r+1}+e_{1,\dots,r+1})$ acts as the null endmorphism on this copy of $\tilde\Delta_r$. It is not
hard to check that the only projections that induce nonzero endomorphisms are  ${1\over 2} (1+e_1\dots
e_{r})$,  ${1\over 2} (1+e_1\dots e_{r})\cdot e_{r+1,r+2}$, ${1\over 2} (1+e_1\dots e_{r})\cdot e_{r+1,r+3}$
and  ${1\over 2} (1+e_1\dots e_{r})\cdot e_{r+2,r+3}$.  Note that the Hermitian product of $\Delta_{r+3}$
restricts to a positive definite inner product on
$\tilde\Delta_{r+3}$ (cf \rf{eq: real inner product}).
Now we will check whether the endomorphisms induced by ${1\over 2} (1+e_1\dots e_{r})$,  ${1\over 2}
(1+e_1\dots e_{r})\cdot e_{r+1,r+2}$, ${1\over 2} (1+e_1\dots e_{r})\cdot e_{r+1,r+3}$ and  ${1\over 2}
(1+e_1\dots e_{r})\cdot e_{r+2,r+3}$ are
symmetric or antisymmetric:
\begin{itemize}
 \item The element ${1\over 2} (1+e_1\dots e_{r})\in Cl_{r+3}$ acts as the identity on this copy of
$\tilde\Delta_{r}$ so is a
symmetric
automorphism.
 \item For $v,w\in\tilde\Delta_{r+3}$, the element ${1\over 2} (1+e_1\dots e_{r})\cdot
e_{r+1,r+2}\in\mathfrak{spin}(r+3)$ is such that
{\footnotesize
\begin{eqnarray*}
   &&
\left< {1\over 2} (1+e_1\dots e_{r})\cdot e_{r+1,r+2}\cdot {1\over 2}
(1+e_1\dots e_r)\cdot(v+\gamma_{r+3}(v)), {1\over 2}
(1+e_1\dots e_r)\cdot(w+\gamma_{r+3}(w)) \right>\\
   &&=
-\left<  {1\over 2}
(1+e_1\dots e_r)\cdot(v+\gamma_{r+3}(v)),
{1\over 2} (1+e_1\dots e_{r})\cdot e_{r+1,r+2}\cdot{1\over 2}
(1+e_1\dots e_r)\cdot(w+\gamma_{r+3}(w)) \right>,
\end{eqnarray*}
}
so that ${1\over 2} (1+e_1\dots e_{r})\cdot e_{r+1,r+2}$ induces a complex structure $I$ on
$\tilde\Delta_{r}$. Indeed, it is a complex
structure.
 \item Similarly, ${1\over 2} (1+e_1\dots e_{r})\cdot e_{r+1,r+3}$ and  ${1\over 2} (1+e_1\dots
e_{r})\cdot e_{r+2,r+3}$ induce complex structures $J$ and $K$ on $\tilde\Delta_{r}$.
\end{itemize}
Thus,
\[C_{\mathfrak{so}(d_r)}(\mathfrak{spin}(r)) = \mathfrak{sp}(1) =\span\left( I, J,
K\right)\subset\mathfrak{so}(d_r)\cong\ext^2\tilde\Delta_r.\]

\subsubsection*{Case $r\equiv -3 \,\,\,({\rm mod}\,\,\, 8)$ }

Recall that
\[\tilde\Delta_r\cong \tilde\Delta_{r+2}\]
as $\mathfrak{spin}(r)$ representations,
and
\[\tilde\Delta_r\otimes \tilde\Delta_r\cong\tilde\Delta_{r+2}\otimes \tilde\Delta_{r+2}.\]
Since $r+2\equiv -1$ (mod 8),
\[\tilde\Delta_{r+2}\otimes \tilde\Delta_{r+2} \cong \ext^{ev}\mathbb{R}^{r+2}\]
as a $\mathfrak{spin}(r+2)$ representation
and $\mathbb{R}^{r+2} = \mathbb{R}^r \oplus 2$,
\begin{eqnarray*}
 \ext^0 \mathbb{R}^{r+2} &=& 1,\\
 \ext^2 \mathbb{R}^{r+2} &=& \ext^2 \mathbb{R}^r + 2 \mathbb{R}^r + 1,\\
 \ext^4 \mathbb{R}^{r+2} &=& \ext^4 \mathbb{R}^r + 2 \ext^3\mathbb{R}^r + \ext^2\mathbb{R}^r,\\
 \vdots\\
 \ext^{r+1} \mathbb{R}^{r+2} &=& \ext^{r+1} \mathbb{R}^r + 2 \ext^r\mathbb{R}^r + \ext^{r-1}\mathbb{R}^r,
\end{eqnarray*}
so that
\begin{eqnarray*}
\tilde\Delta_r\otimes \tilde\Delta_r
  &=&2\ext^*\mathbb{R}^r.
\end{eqnarray*}

Recall that we wish to compute the centralizer of $\mathfrak{spin}(r)$ in 
$\mathfrak{so}(d_r)\subset \End(\tilde\Delta_r)= \End(\tilde\Delta_{r+2})\cong Cl_{r+2}^0$.
By Lemma \ref{Clifford commute}, 
\[C_{Cl_{r+2}^0}(\mathfrak{spin}(r))=\span(1, \quad e_{r+1}e_{r+2},\quad e_1\dots e_{r+1},\quad e_1\dots
e_re_{r+2}), \]
where the last three elements form a copy of $\mathfrak{sp}(1)$.
By means of Clifford multiplication, these three elements act as orthogonal complex structures $I$, $J$, $K$
on $\tilde\Delta_{r}$ and behave as quaternions, i.e. 
\[\span(I,J,K)=\mathfrak{sp}(1)\subset \mathfrak{so}(d_r)\cong
\ext^2\tilde\Delta_r.\]

\subsubsection*{Case $r\equiv 0 \,\,\,({\rm mod}\,\,\, 8)$ }

Recall  that 
$\mathfrak{spin}(r)$  has two irreducible representations given by 
\[\tilde\Delta_r^\pm={1\over 2}(1\pm e_1\dots
e_r)\cdot\tilde\Delta_{r+1},\] 
so that 
\[ \tilde\Delta_{r+1}\cong\tilde\Delta_r^+\oplus \tilde\Delta_r^-,\]
and
\[(\tilde\Delta_r^+\oplus \tilde\Delta_r^-)\otimes (\tilde\Delta_r^+\oplus
\tilde\Delta_r^-)\cong\tilde\Delta_{r+1}\otimes \tilde\Delta_{r+1}\]
as $\mathfrak{spin}(r)$ representations.
Since $r+1\equiv 1$ (mod 8), 
\[\tilde\Delta_{r+1}\otimes \tilde\Delta_{r+1} \cong \ext^{ev}\mathbb{R}^{r+1}\]
as a $\mathfrak{spin}(r+1)$ representation,
and $\mathbb{R}^{r+1} = \mathbb{R}^r \oplus 1$,
\begin{eqnarray*}
 \ext^0 \mathbb{R}^{r+1} &=& 1,\\
 \ext^2 \mathbb{R}^{r+1} &=& \ext^2 \mathbb{R}^r +  \mathbb{R}^r,\\
 \ext^4 \mathbb{R}^{r+1} &=& \ext^4 \mathbb{R}^r + \ext^3\mathbb{R}^r,\\
 \vdots\\
 \ext^{r} \mathbb{R}^{r+1} &=& \ext^{r} \mathbb{R}^r + \ext^{r-1}\mathbb{R}^r,
\end{eqnarray*}
i.e.
\begin{eqnarray*}
(\tilde\Delta_r^+\oplus \tilde\Delta_r^-)\otimes (\tilde\Delta_r^+\oplus \tilde\Delta_r^-)
  &=&\ext^*\mathbb{R}^r,
\end{eqnarray*}
which has only $2$ trivial summands with respect to $\mathfrak{spin}(r)$. On the other hand,
\begin{eqnarray*}
(\tilde\Delta_r^+\oplus \tilde\Delta_r^-)\otimes (\tilde\Delta_r^+\oplus \tilde\Delta_r^-)
  &=& \ext^2\tilde\Delta_r^+ \oplus \sym_0^2 \tilde\Delta_r^+ \oplus 1
\oplus\ext^2\tilde\Delta_r^- \oplus \sym_0^2 \tilde\Delta_r^- \oplus 1
\oplus\tilde\Delta_r^+\otimes\tilde\Delta_r^-\oplus\tilde\Delta_r^-\otimes\tilde\Delta_r^+,
 \end{eqnarray*}
i.e. no other summand contains a trivial $\mathfrak{spin}(r)$ representation.

Recall that we wish to find the centralizer of $\mathfrak{spin}(r)$ in
$\mathfrak{so}(d_r)\oplus \mathfrak{so}(d_r)\subset \End(\tilde\Delta_r^+)\oplus\End(\tilde\Delta_r^-)\cong
Cl_r^0$.
By Lemma \ref{Clifford commute},
\[C_{Cl_{r}^0}(\mathfrak{spin}(r)) = \span(1)\oplus\span(e_1\dots e_r).\]
Since both $1$ and $ e_1\dots e_r$ induce symmetric endomorphisms on $\tilde\Delta_r$,
\[C_{\mathfrak{so}(d_r)\oplus\mathfrak{so}(d_r)}(\mathfrak{spin}(r))=\{0\}\subset\mathfrak{so}
(d_r)\cong\ext^2\tilde\Delta_r.\]

\subsubsection*{Case $r\equiv 4 \,\,\,({\rm mod}\,\,\, 8)$ }

Recall  that  $\mathfrak{spin}(r)$  has two irreducible
representations and
\[\tilde\Delta_r^+\oplus \tilde\Delta_r^-\cong
\tilde\Delta_{r+1}\cong\tilde\Delta_{r+2}\cong\tilde\Delta_{r+3}\]
as representations of $\mathfrak{spin}(r)$.
Since $r+3\equiv -1$ (mod 8), $\gamma_{r+3}$ is a real structure and
\[\tilde\Delta_{r+3}=\{v+\gamma_{r+3}(v) \,\,|\,\, v\in\Delta_{r+3}\}.\]
Moreover 
\[\tilde\Delta_r^\pm={1\over 2}(1\pm e_1\dots e_r)\cdot\tilde\Delta_{r+3}\]
and
\[\tilde\Delta_{r+3}\otimes \tilde\Delta_{r+3}\cong(\tilde\Delta_r^+\oplus \tilde\Delta_r^-)\otimes
(\tilde\Delta_r^+\oplus \tilde\Delta_r^-)\]
with respect to $\mathfrak{spin}(r)$.
With respect to $\mathfrak{spin}(r+3)$,
\[\tilde\Delta_{r+3}\otimes \tilde\Delta_{r+3} \cong \ext^{ev}\mathbb{R}^{r+3}.\]
Now, $\mathbb{R}^{r+3} = \mathbb{R}^r \oplus 3$,
\begin{eqnarray*}
 \ext^0 \mathbb{R}^{r+3} &=& 1,\\
 \ext^2 \mathbb{R}^{r+3} &=& \ext^2 \mathbb{R}^r +  3\mathbb{R}^r+3,\\
 \ext^4 \mathbb{R}^{r+3} &=& \ext^4 \mathbb{R}^r + 3\ext^3\mathbb{R}^r+3\ext^2\mathbb{R}^r+\mathbb{R}^r,\\
 \vdots\\
 \ext^{r+2} \mathbb{R}^{r+3} &=& \ext^{r+2} \mathbb{R}^r + 3\ext^{r+1}\mathbb{R}^r+3\ext^{r}\mathbb{R}^r+
\ext^{r-1}\mathbb{R}^r= 3+ \ext^{r-1}\mathbb{R}^r,
\end{eqnarray*}
and
\begin{eqnarray*}
(\tilde\Delta_r^+\oplus \tilde\Delta_r^-)\otimes (\tilde\Delta_r^+\oplus \tilde\Delta_r^-)
  &=& 4\,\,\ext^*\mathbb{R}^r.
\end{eqnarray*}

Recall that we wish to compute the centralizer of $\mathfrak{spin}(r)$ in 
$\mathfrak{so}(d_r)\oplus\mathfrak{so}(d_r)\subset \End(\tilde\Delta_r^+\oplus
\tilde\Delta_r^-)=\End(\tilde\Delta_{r+3})\cong Cl_{r+3}^0.$
First we will compute
$C_{Cl_{r+3}^0}(\mathfrak{spin}(r))$.
If
\[\eta = \sum_{|I|\equiv 0 (2)} \eta_Ie_I \in C_{Cl_{r+3}^0}(\mathfrak{spin}(r))\]
then it must commute in Clifford product with  every $e_ie_j\in \mathfrak{spin}(r)$, $1\leq i<j\leq r$.
By Lemma \ref{Clifford commute}, 
the only free coefficients are
$\eta_{\emptyset}$,
$\eta_{r+1,r+2}$,
$\eta_{r+1,r+3}$,
$\eta_{r+2,r+3}$,
$\eta_{1,\dots,r}$,
$\eta_{1,\dots,r+2}$,
$\eta_{1,\dots,r+1,r+3}$,
$\eta_{1,\dots,r,r+2,r+3}$,
i.e.
\begin{eqnarray*}
C_{Cl_{r+3}^0}(\mathfrak{spin}(r)) 
&=& \span(1,
e_{r+1,r+2},
e_{r+1,r+3},
e_{r+2,r+3},
e_{1,\dots,r},
e_{1,\dots,r+2},
e_{1,\dots,r+1,r+3},
e_{1,\dots,r,r+2,r+3}
)\\
&=& \span\left(
{1\over 2}(1\pm e_{1\ldots r}),
{1\over 2}(1\pm e_{1\ldots r})e_{r+1}e_{r+2},
{1\over 2}(1\pm e_{1\ldots r})e_{r+1}e_{r+3},
{1\over 2}(1\pm e_{1\ldots r})e_{r+2}e_{r+3}
\right). 
\end{eqnarray*}
Now we need to check which of these elements induce antisymmetric endomorphisms on $\tilde\Delta_r^\pm$.
respectively.
\begin{itemize}
 \item The element ${1\over 2}(1\pm e_{1\ldots r})\in Cl_{r+3}^0$ induces the identity endomorphism on
$\tilde\Delta_r^\pm$ and the null endomorphism on $\tilde\Delta_r^\mp$, both of which are symmetric.
 \item The elements $\frac12(1\pm e_{1}\dots e_{r})e_{r+1,r+2}$, $\frac12(1\pm e_{1}\dots e_{r})e_{r+1,r+3}$
and $\frac12(1\pm e_{1}\dots e_{r})e_{r+2,r+3}$ induce almost complex structures $I^\pm$, $J^\pm$, $K^\pm$ on
$\tilde\Delta_{r}^{\pm}$ respectively, and the null endomorphism on $\tilde\Delta^\mp$.
Such elements also commute with the elements of $\mathfrak{spin}(r)$. In other
words, 
\[\mathfrak{sp}(1)^\pm =  \span(I^\pm,J^\pm,K^\pm)\subset \ext^2\tilde\Delta_{r}^{\pm}\]
are trivial $\mathfrak{spin}(r)$ representations.
\end{itemize}
Hence,
\[C_{\mathfrak{so}(d_r)\oplus\mathfrak{so}(d_r)}(\mathfrak{spin}(r))\cong
\mathfrak{sp}(1)^+\oplus\mathfrak{sp}(1)^-.\]
\qd

\section{Centralizers}\label{sec: centralizers}
Due to geometric considerations in \cite{Moroianu-Semmelmann,Herrera-Santana}, we will consider
$\mathfrak{spin}(r)$ embedded in $\mathfrak{so}(N)$ in the following way. Suppose that $Cl_r^0$ is represented
on $\mathbb{R}^N$, for some $N\in\mathbb{N}$,
in such a way that each bivector $e_ie_j$ is mapped to an antisymmetric endomorphism $J_{ij}$ satisfying
\begin{equation}
J_{ij}^2 = -{\rm Id}_{\mathbb{R}^N}.\label{eq:almost-complex-structures}
\end{equation}

\subsection{Centralizer of $\mathfrak{spin}(r)$ in $\mathfrak{so}(d_rm)$, \,\,$r\not\equiv 0
\,\,\,({\rm mod}\,\,\,\, 4)$, $r>1$} \label{sec: centralizer r not 4}

Let us assume $r\not\equiv 0
\,\,\,({\rm mod}\,\,\,\, 4)$, $r>1$.
In this case, $\mathbb{R}^N$ decomposes into a sum of irreducible representations of $Cl_r^0$.
Since this algebra is simple, such irreducible representations can only be trivial or
copies of the standard representation
$\tilde\Delta_r$ of $Cl_r^0$ (cf. \cite{Lawson}). Due to
\rf{eq:almost-complex-structures}, there are no trivial summands in such a decomposition so that
\begin{eqnarray*}
\mathbb{R}^N
&=& \underbrace{\tilde\Delta_r\oplus\cdots\oplus \tilde\Delta_r}_{m
\,\,\,times} .
\end{eqnarray*}
By restricting to $\mathfrak{spin}(r)\subset Cl_r^0$,
\[\mathbb{R}^N =\tilde\Delta_r \otimes_{\mathbb{R}}\mathbb{R}^m\]
we see that $\mathfrak{spin}(r)$ has an isomorphic image
\[\widehat{\mathfrak{spin}(r)}=\mathfrak{spin}(r)\otimes \{{\rm Id}_{m\times m}\}\subset
\mathfrak{so}(d_rm),\]
which is the subalgebra of $\mathfrak{so}(d_rm)$ whose centralizer 
$C_{\mathfrak{so}(d_rm)}(\widehat{\mathfrak{spin}({r})})$ we wish to find.

\begin{theo}\label{theo: 1}Let $r\not\equiv0$ $(\mbox{\rm{mod} } 4) $ and let
$\widehat{\mathfrak{spin}(r)}
\subset
\mathfrak{so}(d_rm)$ as described before. The centralizer of $\widehat{\mathfrak{spin}(r)}$ in
$\mathfrak{so}(d_rm)$ is isomorphic to
\[  
\begin{array}{|c|c|}
 \hline
  r \mbox{\ {\rm (mod 8)} }  &C_{\mathfrak{so}(d_rm)}(\widehat{\mathfrak{spin}({r})}) \rule{0pt}{3ex}\\
 \hline
 1 & \mathfrak{so}(m)\\
\hline
 2 & \mathfrak{u}(m)\\
\hline
 3 & \mathfrak{sp}(m)\\
\hline
 5 & \mathfrak{sp}(m)\\
\hline
 6 & \mathfrak{u}(m)\\
\hline
 7 & \mathfrak{so}(m)\\
\hline
  \end{array}
\]
\end{theo}
{\em Proof}. Consider the real $(d_rm)$-dimensional real Grassmannian
\[\mathcal{G}={SO(d_r+m)\over SO(d_r)\times SO(m)}.\]
The tangent space factors as follows
\begin{eqnarray*}
T_{[{\rm Id}_{(d_r+m)\times(d_r+m)}]}\mathcal{G}&\cong&\mathbb{R}^{d_r}\otimes\mathbb{R}^m\\
&\cong&\mathbb{R}^{d_rm}.
\end{eqnarray*}
so that the differential of the isotropy representation is
\begin{eqnarray*}
\mathfrak{so}(d_r)\oplus \mathfrak{so}(m) &\longrightarrow&
[\mathfrak{so}(d_r)\otimes \{{\rm Id}_{m\times m}\}]
\oplus
[\{{\rm Id}_{d_r\times d_r}\}\otimes\mathfrak{so}(m)]
\subset
\mathfrak{so}(d_rm)\\
 (A,B) &\mapsto& A\otimes {\rm Id}_{m\times m} \oplus {\rm Id}_{d_r\times d_r}\otimes B.
\end{eqnarray*}
Let $\widehat{\mathfrak{so}(m)}=\{{\rm Id}_{d_r\times d_r}\}\otimes\mathfrak{so}(m)$
and
$\widehat{\mathfrak{so}(d_r)}=\mathfrak{so}(d_r)\otimes\{{\rm Id}_{m\times m}\}$.
Thus, we see that $\widehat{\mathfrak{so}(m)}$ centralizes $\widehat{\mathfrak{so}(d_r)}$ in
$\mathfrak{so}(d_rm)$, and
\[\widehat{\mathfrak{so}(m)}\subseteq C_{\mathfrak{so}(d_rm)}(\widehat{\mathfrak{spin}(r)}).\]
Let us consider the following orthogonal decomposition
\[\mathfrak{so}(d_rm)= [\widehat{\mathfrak{so}(m)}\oplus \widehat{\mathfrak{so}(d_r)}] \oplus \mathfrak{m},\]
and set
\begin{eqnarray*}
 \mathfrak{g} &=& \mathfrak{so}(d_rm),\\
 \mathfrak{h} &=& \widehat{\mathfrak{so}(m)}\oplus\widehat{\mathfrak{so}(d_r)}.
\end{eqnarray*}
Since the homogeneous space
\[\mathcal{F}={SO(d_rm)\over SO(d_r)\otimes SO(m)}\]
is Riemannian homogeneous, it is reductive, i.e.
\begin{eqnarray*}
{} [\mathfrak{h},\mathfrak{m}] &\subset& \mathfrak{m}.
\end{eqnarray*}
Let
\[X=X_1+X_2+X_3 \in \mathfrak{g}\]
where
\begin{eqnarray*}
 X_1 &\in&  \widehat{\mathfrak{so}(m)} ,\\
 X_2 &\in&  \widehat{\mathfrak{so}(d_r)} ,\\
 X_3 &\in&  \mathfrak{m},
\end{eqnarray*}
and assume that $X\in C_{\mathfrak{so}(d_rm)}(\widehat{\mathfrak{spin}(r)})$, i.e.
\[[X,Y]=0\]
for all $Y\in \widehat{\mathfrak{spin}(r)}$.
Thus,
\begin{eqnarray*}
 0
  &=& [X_1,Y]+[X_2,Y]+[X_3,Y].
\end{eqnarray*}
Note that
\begin{eqnarray*}
{} [X_1,Y] &\in& \mathfrak{h},\\
{} [X_2,Y] &\in& \mathfrak{h},\\
{} [X_3,Y] &\in& \mathfrak{m},
\end{eqnarray*}
so that
\begin{eqnarray*}
{} [X_1+X_2,Y] &=&0,\\
{} [X_3,Y] &=&0.
\end{eqnarray*}
Since $X_1\in\widehat{\mathfrak{so}(m)}$ and
$Y\in\widehat{\mathfrak{spin}(r)}\subset\widehat{\mathfrak{so}(d_r)}$,
\[[X_1,Y]=0,\]
which implies
\[[X_2,Y]=0.\]
On the other hand, since
\[[X_3,Y]=0\]
for all $Y\in \widehat{\mathfrak{spin}(r)}$,
the subalgebra $\widehat{\mathfrak{spin}(r)}\subset\mathfrak{h}$ acts trivially on the 1-dimensional subspace
of the tangent space $\mathfrak{m}$ of $\mathcal{F}$ at $[{\rm Id}_{(d_rm)\times (d_rm)}]$
generated by $X_3$. Now, as a representation of
$\mathfrak{h}=\widehat{\mathfrak{so}(d_r)}\oplus \widehat{\mathfrak{so}(m)}\cong \mathfrak{so}(d_r)\oplus
\mathfrak{so}(m)$,
\[\mathfrak{m} \cong \left[\ext^2 \mathbb{R}^{d_r}\otimes \sym_0^2\mathbb{R}^{m}\right]
\oplus
\left[\sym_0^2\mathbb{R}^{d_r}\otimes\ext^2 \mathbb{R}^m\right] .\]
By restricting to $\widehat{\mathfrak{so}(d_r)}$
\begin{eqnarray*}
\mathfrak{m}
&\cong&
\left[\ext^2 \mathbb{R}^{d_r}\otimes \left({m+1\choose2}-1\right)\right]\oplus
\left[\sym_0^2\mathbb{R}^{d_r}\otimes{m\choose2}\right] ,
\end{eqnarray*}
i.e. $\mathfrak{m}$ decomposes as the sum of multiple copies of the
irreducible
$\mathfrak{so}(d_r)$
representations $\ext^2\mathbb{R}^{d_r}$ and
$\sym_0^2\mathbb{R}^{d_r}$.
By restricting further to
$\widehat{\mathfrak{spin}(r)}\subset\widehat{\mathfrak{so}(d_r)}$,
$\mathfrak{m}$ decomposes as

\begin{eqnarray}
\mathfrak{m}
&\cong&
\left[\ext^2 \tilde\Delta_r\otimes \left({m+1\choose2}-1\right)\right]\oplus
\left[\sym_0^2\tilde\Delta_r\otimes{m\choose2}\right] .\label{eq: decomposition of m v1}
\end{eqnarray}
Both $\mathfrak{spin}(r)$ representations $\ext^2 \tilde\Delta_r$ 
and $\sym_0^2\tilde\Delta_r$ decompose further into irreducible summands. 
Now we need to work out three cases separately.

\subsubsection*{Case $r\equiv\pm1 \mbox{ {\rm (mod 8)}}$:}
By Proposition \ref{prop: 1}, the centralizer of
$\widehat{\mathfrak{spin}(r)}$ in $\widehat{\mathfrak{so}(d_r)}$ is trivial, i.e.
\[X_2=0.\]
Recall that $\widehat{\mathfrak{spin}(r)}$ preserves each summand in \rf{eq: decomposition of m v1} and
annihilates $X_3$.  
By Proposition \ref{prop: 1}, 
there are
no trivial summands in either $\ext^2\tilde\Delta_r$ nor $\sym_0^2\tilde\Delta_r$, i.e.
\[X_3=0.\]
Hence
\[X=X_1\in\widehat{\mathfrak{so}(m)}.\]

\subsubsection*{Case $r\equiv\pm2 \mbox{ {\rm (mod 8)}}$:}
By Proposition \ref{prop: 1}, the centralizer of
$\widehat{\mathfrak{spin}(r)}$ in $\widehat{\mathfrak{so}(d_r)}$ is a copy of $\mathfrak{u}(1)$, i.e.
\[X_2=\lambda \,J\otimes {\rm Id}_{m\times m},\]
where $J$ is an orthogonal complex structure that generates $\mathfrak{u}(1)$ and  $\lambda \in \mathbb{R}$.
Recall that $\widehat{\mathfrak{spin}(r)}$ preserves each summand in \rf{eq: decomposition of m v1} and
annihilates $X_3$.  
There are no trivial summands in $\sym_0^2\tilde\Delta_r$,
but there is a trivial summand in $\ext^2\tilde\Delta_r$
generated precisely by
$J$, since it is an antisymmetric endomorphism.
We see that $\mathfrak{m}$ contains
\[\span(J)\otimes\sym_0^2\mathbb{R}^{m}\]
as a trivial $\widehat{\mathfrak{spin}(r)}$ representation.
Hence
\[X\in  \widehat{\mathfrak{so}(m)} \oplus \span(J)\otimes(\span({\rm Id}_{m\times m})\oplus
\sym_0^2\mathbb{R}^m) \subset
\mathfrak{so}(d_rm).\]

In order to recogize which Lie algebra $\left[\widehat{\mathfrak{so}(m)} \oplus \span(J)\otimes
\sym^2\mathbb{R}^m\right]$ is, notice that if $A\in \mathfrak{u}(m)$, by separating real and imaginary parts
\[A= A_1 + i A_2,\]
$A_1\in \mathfrak{so}(m)$ is antisymmetric and
$A_2$ is symmetric, i.e. $A_2\in\sym^2 \mathbb{R}^m$.
Here, a canonical summand $\mathfrak{u}(1)$ is spanned by the element $i{\rm Id}_{m\times m}$. Note that due
to the existence of $J$, we can work instead with a complex vector space, where
$J$ corresponds to $i$, $J\otimes \sym_0^2\mathbb{R}^m$ corresponds to $i\sym_0^2\mathbb{R}^m$ and
\[ \widehat{\mathfrak{so}(m)} \oplus \span(J)\otimes
\sym^2\mathbb{R}^m \cong \mathfrak{u}(m). \]

\subsubsection*{Case $r\equiv\pm3 \mbox{ {\rm (mod 8)}}$:}

By Proposition \ref{prop: 1}, the centralizer of $\widehat{\mathfrak{spin}(r)}$ in
$\widehat{\mathfrak{so}(d_r)}$ is a copy of $\mathfrak{sp}(1)=\span(I,J,K)$,
where $I,J,K$ are three orthogonal complex structures which behave as imaginary quaternions.
Thus, 
\[X_2\in\mathfrak{sp}(1)\otimes \span({\rm Id}_{m\times m})
.\]
By Proposition \ref{prop: 1}, $\sym_0^2\tilde\Delta_{r}$ contains no trivial $\widehat{\mathfrak{spin}(r)}$
representations, but 
$\ext^2\tilde\Delta_{r}$ does contain a $3$-dimensional one given by $\mathfrak{sp}(1)=\span(I,J,K)$.
We have the trivial $\widehat{\mathfrak{spin}(r)}$ representation in $\mathfrak{m}$
\[
\span(I,J,K)\otimes\sym_0^2\mathbb{R}^{m}
=\mathfrak{sp}(1)\otimes\sym_0^2\mathbb{R}^{m}.
\]
Altogether, we have that
\[X\in 
=
\mathfrak{so}(m) \oplus \mathfrak{sp}(1)\otimes \sym^2\mathbb{R}^{m}.
\]
In order to recognize this Lie algebra, notice that if $A\in \mathfrak{sp}(m)$, by separating real and
imaginary parts
\[A= A_1 + i A_2 + j A_3+ k A_4,\]
$A_1\in \mathfrak{so}(m)$ is antisymmetric and
$A_2, A_3, A_4$ are symmetric, i.e. $A_2, A_3, A_4\in\sym^2 \mathbb{R}^m$.
The summand $\mathfrak{sp}(1)$ is spanned by the elements    $i{\rm Id}_{m\times m}, j{\rm Id}_{m\times m}, k{\rm
Id}_{m\times m}$.
Moreover,
due to the existence of $I,J,K$, we can work instead with a quaternionic vector space, in which, $I$
corresponds to $i$, $J$ corresponds to $j$, $K$ corresponds to $k$, and
$\span(I,J,K)\otimes \sym_0^2\mathbb{R}^m$ corresponds to
$i\sym_0^2\mathbb{R}^m\oplus j\sym_0^2\mathbb{R}^m\oplus k\sym_0^2\mathbb{R}^m$ so that
\[ \widehat{\mathfrak{so}(m)} \oplus \mathfrak{sp}(1)\otimes
\sym^2\mathbb{R}^m \cong \mathfrak{sp}(m). \]
\qd

\subsection{Centralizer of $\mathfrak{spin}(r)$ in $\mathfrak{so}(d_rm_1+d_rm_2)$, $r\equiv 0
\,\,\,({\rm mod}\,\,\,\, 4)$}\label{sec:centralizer-Spin(4k)}

Let us assume $r\equiv 0
\,\,\,({\rm mod}\,\,\,\, 4)$.
Recall that if $\hat\Delta_r$ is the irreducible representation of $Cl_r$, then by restricting this
representation to $Cl^0_r$ it splits as the sum of two inequivalent irreducible representations
\[\hat\Delta_r = \tilde{\Delta}_r^+ \oplus \tilde{\Delta}_r^-.\]
Since $\mathbb{R}^N$ is a representation of $Cl_r^0$ satisfying
\rf{eq:almost-complex-structures}, there are no trivial summands in such a decomposition so that
\begin{eqnarray*}
\mathbb{R}^N
&=& \tilde\Delta_r^+\otimes \mathbb{R}^{m_1} \oplus \tilde\Delta_r^-\otimes \mathbb{R}^{m_2}.
\end{eqnarray*}
By restricting this representation to $\mathfrak{spin}(r)\subset Cl_r^0$,
consider
\[\widehat{\mathfrak{spin}(r)}=\mathfrak{spin}(r)^+\otimes ({\rm Id}_{m_1\times m_1}\oplus
\mathbf{0}_{m_2\times m_2})\oplus
\mathfrak{spin}(r)^-\otimes(\mathbf{0}_{m_1\times m_1}\oplus {\rm Id}_{m_2\times m_2})
\subset
\mathfrak{so}(d_rm_1+d_rm_2),\]
where $\mathfrak{spin}(r)^{\pm}$ are the images of $\mathfrak{spin}(r)$ in 
$\End(\tilde\Delta_r^{\pm})$
respectively.
We wish to find the centralizer 
$C_{\mathfrak{so}(d_rm_1+d_rm_2)}(\widehat{\mathfrak{spin}({r})})$.

\begin{theo}\label{theo: 2} 
Let $r\equiv 0 \mbox{ {\em (mod 4)} }$. The centralizer of
$\widehat{\mathfrak{spin}(r)}$ in
$\mathfrak{so}(d_rm_1+d_rm_2)$ is isomorphic to
\[ 
\begin{array}{|c|c|}
 \hline
  r \mbox{ {\rm (mod 8)} }   &C_{\mathfrak{so}(d_rm_1+d_rm_2)}(\widehat{\mathfrak{spin}({r})})
\rule{0pt}{3ex}\\
 \hline
    0 &\mathfrak{so}(m_1)\oplus\mathfrak{so}(m_2)\\
\hline
    4 & \mathfrak{sp}(m_1)\oplus\mathfrak{sp}(m_2)\\
\hline
  \end{array}
\]
\end{theo}
{\em Proof}. Consider the homogeneous space
\[\mathcal{G}={SO(m_1+d_r)\times SO(m_2+d_r)\over (SO(d_r)\times SO(m_1))\times (SO(d_r)\times SO(m_2)) }.\]
with the obvious inclusions of subgroups. The tangent space decomposes as follows
\begin{eqnarray*}
T_{[{\rm Id}_{(2d_r+m_1+m_2)\times(d_r+m_1+m_2)}]}\mathcal{G}&\cong&
\mathbb{R}^{m_1}\otimes\mathbb{R}^{d_r}\oplus \mathbb{R}^{m_2}\otimes\mathbb{R}^{d_r},
\end{eqnarray*}
Let 
\begin{eqnarray*}
\widehat{\mathfrak{so}(m_1)}&=&({\rm Id}_{d_r\times d_r}\oplus \mathbf{0}_{d_r\times
d_r})\otimes\mathfrak{so}(m_1) ,\\
\widehat{\mathfrak{so}(m_2)}&=&(\mathbf{0}_{d_r\times d_r}\oplus{\rm Id}_{d_r\times
d_r})\otimes\mathfrak{so}(m_2),\\
\widehat{\mathfrak{so}(d_r)_1}&=&\mathfrak{so}(d_r)\otimes({\rm Id}_{m_1\times
m_1}\oplus \mathbf{0}_{m_2\times m_2}),\\
\widehat{\mathfrak{so}(d_r)_2}&=&\mathfrak{so}(dr)\otimes(\mathbf{0}_{m_1\times m_1}\oplus {\rm
Id}_{m_2\times m_2}).
\end{eqnarray*}
We see that $\widehat{\mathfrak{so}(m_1)}\oplus\widehat{\mathfrak{so}(m_2)}$
centralizes
$\widehat{\mathfrak{so}(d_r)_1}\oplus \widehat{\mathfrak{so}(d_r)_2}$ in
$\mathfrak{so}(d_rm_1+d_rm_2)$,
\[\widehat{\mathfrak{so}(m_1)}\oplus\widehat{\mathfrak{so}(m_2)}\subseteq
C_{\mathfrak{so}(d_rm_1+d_rm_2)}(\widehat{\mathfrak{spin}(r)}).\]
Let us consider the following orthogonal decomposition
\[\mathfrak{so}(d_rm_1+d_rm_2)= [\widehat{\mathfrak{so}(m_1)}\oplus \widehat{\mathfrak{so}(d_r)}_1] \oplus
[\widehat{\mathfrak{so}(m_2)}\oplus \widehat{\mathfrak{so}(d_r)_2}] \oplus
\mathfrak{m},\]
and set
\begin{eqnarray*}
 \mathfrak{g} &=& \mathfrak{so}(d_rm_1+d_rm_2),\\
 \mathfrak{h} &=&
\widehat{\mathfrak{so}(m_1)}\oplus\widehat{\mathfrak{so}(d_r)_1} \oplus
\widehat{\mathfrak{so}(m_2)}\oplus\widehat{\mathfrak{so}(d_r)_2}.
\end{eqnarray*}
Since the homogeneous space
\[{\mathcal F}={SO(d_rm_1+d_rm_2)\over (SO(m_1)\otimes SO(d_r))\times (SO(m_2)\otimes SO(d_r))}\]
is Riemannian homogeneous, it is reductive, and
\begin{eqnarray*}
{} [\mathfrak{h},\mathfrak{m}] &\subset& \mathfrak{m}.
\end{eqnarray*}
Let
\[X=X_1+X_2+X_3 \in \mathfrak{g}\]
where
\begin{eqnarray*}
 X_1 &\in&  \widehat{\mathfrak{so}(m_1)}\oplus\widehat{\mathfrak{so}(m_2)} ,\\
 X_2 &\in&  \widehat{\mathfrak{so}(d_r)_1} \oplus\widehat{\mathfrak{so}(d_r)_2} ,\\
 X_3 &\in&  \mathfrak{m},
\end{eqnarray*}
and assume that $X\in C_{\mathfrak{so}(d_rm_1+d_rm_2)}(\widehat{\mathfrak{spin}({r})})$, i.e.
\[[X,Y]=0\]
for all $Y\in \widehat{\mathfrak{spin}(r)}$. Thus,
\begin{eqnarray*}
 0
  &=& [X_1,Y]+[X_2,Y]+[X_3,Y].
\end{eqnarray*}
Note that
\begin{eqnarray*}
{} [X_1,Y] &\in& \mathfrak{h},\\
{} [X_2,Y] &\in& \mathfrak{h},\\
{} [X_3,Y] &\in& \mathfrak{m},
\end{eqnarray*}
so that
\begin{eqnarray*}
{} [X_1+X_2,Y] &=&0,\\
{} [X_3,Y] &=&0.
\end{eqnarray*}
Since $X_1\in\widehat{\mathfrak{so}(m_1)}\oplus\widehat{\mathfrak{so}(m_2)}$ and
$Y\in\widehat{\mathfrak{spin}(r)}\subset\widehat{\mathfrak{so}(d_r)_1}\oplus\widehat{\mathfrak{so}(d_r)_2}$,
\[[X_1,Y]=0,\]
which implies
\[[X_2,Y]=0.\]
Since
\[[X_3,Y]=0\]
for all $Y\in \widehat{\mathfrak{spin}(r)}$,
the subalgebra $\widehat{\mathfrak{spin}(r)}\subset\mathfrak{h}$ acts trivially on the 1-dimensional subspace
of the tangent space $\mathfrak{m}$ of $\mathcal{F}$ at $[{\rm Id}_{(d_rm_1+d_rm_2)\times(d_rm_1+d_rm_2)}]$
generated by $X_3$.
Note that
\begin{eqnarray*}
 \mathfrak{so}(d_rm_1+d_rm_2)
   &=&
  \ext^2 (\mathbb{R}_1^{d_r}\otimes \mathbb{R}^{m_1} \oplus \mathbb{R}_2^{d_r}\otimes \mathbb{R}^{m_2})\\
   &=&
  \ext^2 (\mathbb{R}_1^{d_r}\otimes \mathbb{R}^{m_1})
\oplus   (\mathbb{R}_1^{d_r}\otimes \mathbb{R}^{m_1}) \otimes (\mathbb{R}_2^{d_r}\otimes \mathbb{R}^{m_2})
 \oplus \ext^2  (\mathbb{R}_2^{d_r}\otimes \mathbb{R}^{m_2})\\
   &\cong&
  \widehat{\mathfrak{so}(d_r)_1} \oplus \widehat{\mathfrak{so}(m_1)} \oplus
\left[\ext^2 \mathbb{R}_1^{d_r}\otimes \sym_0^2\mathbb{R}^{m_1} \oplus
\sym_0^2 \mathbb{R}_1^{d_r}\otimes\ext^2 \mathbb{R}^{m_1}\right] \\
   &&
\oplus\,\,
 \mathbb{R}_1^{d_r}\otimes \mathbb{R}_2^{d_r}\otimes \mathbb{R}^{m_1}\otimes \mathbb{R}^{m_2} \\
   &&
  \oplus\,\,
\widehat{\mathfrak{so}(d_r)_2} \oplus \widehat{\mathfrak{so}(m_2)} \oplus
\left[\ext^2 \mathbb{R}_2^{d_r}\otimes \sym_0^2\mathbb{R}^{m_2} \oplus
\sym_0^2 \mathbb{R}_2^{d_r}\otimes\ext^2 \mathbb{R}^{m_2}\right] ,
\end{eqnarray*}
so that, by restricting to $\widehat{\mathfrak{spin}(r)}$,
\begin{eqnarray*}
 \mathfrak{m}
   &=&
\left[\ext^2 \tilde\Delta_r^+\otimes \sym_0^2\mathbb{R}^{m_1} \oplus
\sym_0^2  \tilde\Delta_r^+\otimes\ext^2 \mathbb{R}^{m_1}\right] \\
   &&
\oplus\,\,
 \tilde\Delta_r^+\otimes \tilde\Delta_r^-\otimes \mathbb{R}^{m_1}\otimes \mathbb{R}^{m_2} \\
   &&
  \oplus\,\,
\left[\ext^2 \tilde\Delta_r^-\otimes \sym_0^2\mathbb{R}^{m_2} \oplus
\sym_0^2  \tilde\Delta_r^-\otimes\ext^2 \mathbb{R}^{m_2}\right].
\end{eqnarray*}
Now we need to check two cases separately.

\subsubsection*{Case $r\equiv0 \mbox{ {\rm (mod 8)}}$:}

By Proposition \ref{prop: 1}, the centralizer of $\widehat{\mathfrak{spin}(r)}$ in
$\widehat{\mathfrak{so}(d_r)_1}\oplus\widehat{\mathfrak{so}(d_r)_2}$ is
trivial, i.e. 
\[X_2=0.\]
By Proposition \ref{prop: 1}, $\frak{m}$ has no trivial summands, i.e. 
\[X_3=0.\]
Hence 
\[X=X_1\in \widehat{\mathfrak{so}(m_1)}\oplus\widehat{\mathfrak{so}(m_2)}.\]

\subsubsection*{Case $r\equiv 4 \mbox{ {\rm (mod 8)}}$:}

By Proposition \ref{prop: 1}, the centralizer of $\widehat{\mathfrak{spin}(r)}$ in
$\widehat{\mathfrak{so}(d_r)_1}\oplus\widehat{\mathfrak{so}(d_r)_2}$ is a copy of
$\mathfrak{sp}(1)\oplus\mathfrak{sp}(1)$, i.e.
\[X_2\in\left[\mathfrak{sp}(1)\otimes({\rm Id}_{m_1\times
m_1}\oplus \mathbf{0}_{m_2\times m_2})\right]\oplus\left[\mathfrak{sp}(1)\otimes(\mathbf{0}_{m_1\times
m_1}\oplus{\rm Id}_{m_2\times m_2} )\right]\]
By Proposition \ref{prop: 1},
the only $\widehat{\mathfrak{spin}(r)}$ representations in $\mathfrak m$ containing trivial
$\mathfrak{spin}(r)$ summands are
$\ext^2\tilde\Delta_{r}^\pm$. More precisely, $\ext^2\tilde\Delta_r^\pm$ 
constains a 3-dimensional trivial $\mathfrak{spin}(r)$ representation $\mathfrak{sp}(1)^\pm=\span(I^\pm, 
J^\pm, K^\pm)$, where $I^\pm, J^\pm, K^\pm $ are orthogonal complex structures on $\tilde\Delta_r^\pm$ which
behave as quaternions.
Thus,
we have the trivial $\mathfrak{spin}(r)$ representation in $\mathfrak{m}$
\[\mathfrak{sp}(1)^+\otimes \sym_0^2 \mathbb{R}^{m_1}\oplus\mathfrak{sp}(1)^-\otimes \sym_0^2
\mathbb{R}^{m_2}.\]
Altogether, we have 
\[X\in[\widehat{\mathfrak{so}(m_1)}\oplus\mathfrak{sp}(1)^+\otimes \sym^2
\mathbb{R}^{m_1}]\oplus[\widehat{\mathfrak{so}(m_2)}\oplus\mathfrak{sp}(1)^-\otimes \sym^2
\mathbb{R}^{m_2}]\cong \mathfrak{sp}(m_1)\oplus \mathfrak{sp}(m_2).\]
\qd

{\small
\renewcommand{\baselinestretch}{0.5}
\newcommand{\bi}{\vspace{-.05in}\bibitem} }

\end{document}